# DETECTION OF EDGES IN SPECTRAL DATA II. NONLINEAR ENHANCEMENT[*]

ANNE GELB[†] AND EITAN TADMOR[‡]

**Abstract.** We discuss a general framework for recovering edges in piecewise smooth functions with finitely many jump discontinuities, where $[f](x) := f(x+) - f(x-) \neq 0$. Our approach is based on two main aspects—*localization* using appropriate concentration kernels and *separation of scales* by nonlinear enhancement.

To detect such edges, one employs concentration kernels, $K_\epsilon(\cdot)$, depending on the small scale $\epsilon$. It is shown that *odd kernels, properly scaled, and admissible* (in the sense of having small $W^{-1,\infty}$-moments of order $\mathcal{O}(\epsilon)$) satisfy $K_\epsilon * f(x) = [f](x) + \mathcal{O}(\epsilon)$, thus recovering both the location and amplitudes of all edges. As an example we consider general concentration kernels of the form $K_N^\sigma(t) = \sum \sigma(k/N) \sin kt$ to detect edges from the first $1/\epsilon = N$ spectral modes of piecewise smooth $f$'s. Here we improve in generality and simplicity over our previous study in [A. Gelb and E. Tadmor, *Appl. Comput. Harmon. Anal.*, 7 (1999), pp. 101–135]. Both periodic and nonperiodic spectral projections are considered. We identify, in particular, a new family of exponential factors, $\sigma^{exp}(\cdot)$, with superior localization properties.

The other aspect of our edge detection involves a nonlinear enhancement procedure which is based on separation of scales between the edges, where $K_\epsilon * f(x) \sim [f](x) \neq 0$, and the smooth regions where $K_\epsilon * f = \mathcal{O}(\epsilon) \sim 0$. Numerical examples demonstrate that by coupling concentration kernels with nonlinear enhancement one arrives at effective edge detectors.

## 1. Introduction.

We discuss a general framework for recovering edges from the spectral projections of piecewise smooth functions. Our approach for edge detection is based on two fundamental aspects—localization to the neighborhood of the edges using appropriate concentration kernels and separation of scales by nonlinear enhancement. Both the locations and amplitudes of all edges are recovered.

Let $S_N f(x)$ denote the spectral projection of a piecewise smooth $f$. Given $S_N f$, one can accurately reconstruct $f$ away from its discontinuous jumps, e.g., [10], [14, sect. 2.1], as well as up to the discontinuities [11]. In either case, an a priori knowledge of the location of the edges and their amplitudes is required. This issue was treated in recent literature; consult [1], [5], [13], [15]. In [7], we unified several types of treatments as special cases of appropriate concentration kernels, specifically those discussed in [1] and [13]. Here we improve on these results in both generality and simplicity. To this end, let $[f](x) := f(x+) - f(x-)$ denote the local jump function and let us consider a concentration kernel $K_\epsilon(\cdot)$, depending on a small scale $\epsilon$. It is shown that *odd kernels, properly scaled, and admissible* (in the sense of having small $W^{-1,\infty}$-moments of order $\mathcal{O}(\epsilon)$; see eq. (2.5)), recover both the locations and the amplitudes of the jumps so that

$$(1.1) \qquad\qquad K_\epsilon * f(x) = [f](x) + \mathcal{O}(\epsilon).$$

Thus, $K_\epsilon$ tends to "concentrate" near the singular support of $f$.

---

[*]Received by the editors July 19, 1999; accepted for publication (in revised form) May 16, 2000.

[†]Department of Mathematics, Arizona State University, P.O. Box 871804, Tempe, AZ 85287-1804 (ag@math.la.asu.edu). This author's research was supported in part by the Sloan Foundation.

[‡]Department of Mathematics, UCLA, Los Angeles, CA 90095 (tadmor@math.ucla.edu). This author's research was supported by NSF grant DMS97-06827 and ONR grant N00014-1-J-1076.

Differentiation of $\epsilon$-supported mollifiers is one example for *local* concentration kernels outlined in section 2.2.1. In section 2 we also address the issue of detecting edges in *global* Fourier projections. Given the first $1/\epsilon = N$ Fourier modes, we seek concentration kernels of the form

$$K_N^\sigma(t) = -\sum_{k=1}^{N} \sigma\left(\frac{k}{N}\right)\sin kt.$$

It is shown that if the concentration factors $\sigma(\xi) \equiv \xi\mu(\xi)$ are normalized so that $\int \mu(\xi)d\xi = 1$, then $K_N^\sigma(t)$ is an admissible concentration kernel, $K_N^\sigma * S_N f(x) \to [f](x)$, and the following error estimate holds:

$$\left|\frac{\pi}{N}\sum_{k=-N}^{N}\mu\left(\frac{|k|}{N}\right)\hat{f}_k ike^{ikx} - [f](x)\right| \leq Const \cdot \frac{\log N}{N}.$$

The nonperiodic case is studied in section 3. The analogous results for the Chebyshev case (consult Corollary 3.2) are written as

$$\left|\frac{\pi\sqrt{1-x^2}}{N}\sum_{k=1}^{N}\mu\left(\frac{k}{N}\right)\hat{f}_k T_k'(x) - [f](x)\right| \leq Const \cdot \frac{\log N}{N}.$$

The special cases of Fourier concentration factors $\sigma_\alpha(\xi) \sim \sin\alpha\xi$ and $\sigma^p(\xi) = p\xi^p$ were considered earlier in [1], [7], [9], [13], and [15]. Our general framework motivates a new set of $C_0^\infty$-exponential concentration factors which yield superior localization properties away from the detected edges.

While (1.1) refers to the *asymptotic behavior* of the concentration kernel as a function of the small parameter $\epsilon \downarrow 0$, it is essential to recover the *exact* locations of the edges of $f$ for the accurate reconstruction of $f$. In section 4 we discuss another essential aspect of edge detection, namely, *nonlinear enhancement*. To this end, one introduces a critical threshold, $J_{crit}$, for the amplitude of admissible edges and an enhancement exponent, $p$, to amplify the separation of scales in (1.1) between the edges, where $K_\epsilon * f(x) \sim [f](x) \neq 0$, and the smooth regions where $K_\epsilon * f(x) = \mathcal{O}(\epsilon) \sim 0$. Consider the enhanced kernel

$$K_{\epsilon,J}[f](x) = \begin{cases} K_\epsilon * f & \text{if } \epsilon^{-p/2}|K_\epsilon * f(x)|^p > J_{crit}, \\ 0 & \text{otherwise.} \end{cases}$$

Clearly, with $p$ large enough, one ends up with a sharp edge detector where $K_{\epsilon,J}[f](x) = 0$ at all but $\mathcal{O}(\epsilon)$-neighborhoods of the jump discontinuities. In this sense, the enhancement procedure actually "pinpoints" the location jump discontinuities, allowing an accurate reconstruction of $f$. The particular case $p = 2$ corresponds to the quadratic filter studied in [12], [22], in the special context of concentration kernels based on localized mollifiers.

## 2. Edge detection by concentration kernels.

### 2.1. Concentration kernels.
We want to detect the edges in piecewise smooth functions. Assume that $f(\cdot)$ has jump discontinuities of the first kind with well-defined one-sided limits, $f(x\pm) = \lim_{x\to x\pm} f(x)$, and let $[f](x) := f(x+) - f(x-)$ denote the local jump function. By *piecewise smoothness* we mean[1]

$$(2.1) \qquad F_x(t) := \frac{f(x+t) - f(x-t) - [f](x)}{t} \in BV[0,\delta] \qquad \forall \, x.$$

In practice one encounters functions $f(x)$ with finitely many jump discontinuities, and (2.1) requires the differential of $f(x)$ on each side of the discontinuity to have bounded variation. For example, if $f'(x\pm)$ are well defined (for finitely many jumps), then (2.1) holds.

We will detect the edges in such piecewise smooth $f$'s using smooth *concentration kernels*, $K(t) = K_\epsilon(t)$, depending on a small parameter $\epsilon$. Such kernels are characterized by

$$(2.2) \qquad K_\epsilon * f(x) \longrightarrow [f](x) \qquad \text{as } \epsilon \to 0.$$

Thus the support of $K_\epsilon * f(x)$ tends to "concentrate" near the edges of $f(x)$. One recovers both the location of the jump discontinuities as well as their amplitudes.

To guarantee the concentration property of $K_\epsilon$, we seek *odd* kernels,

$$(2.3) \qquad K_\epsilon(-t) = -K_\epsilon(t),$$

which are normalized so that

$$(2.4) \qquad \int_{t \geq 0} K_\epsilon(t)dt = -1 + \mathcal{O}(\epsilon)$$

and which satisfy the main *admissibility requirement*

$$(2.5) \qquad \left| \int_t t K_\epsilon(t)\phi(t)dt \right| \leq Const \cdot \epsilon ||\phi||_{BV}.$$

*Remarks.*

1. For example, if $K_\epsilon(t)$ concentrates near the origin so that its first moment does not exceed

$$(2.6) \qquad \int_t |t K_\epsilon(t)|dt \leq Const \cdot \epsilon,$$

then it is clearly admissible in the sense that (2.5) holds. We note that our admissibility condition also allows for more general *oscillatory* kernels, $K_\epsilon(t)$, where (2.6) might fail, yet (2.5) is satisfied due to the cancellation effect of the oscillations. Consult (2.17) below.

2. Observe that the admissibility requirement (2.5) generalizes both properties $\mathcal{P}_3$ and $\mathcal{P}_4$ in the definition of admissible kernel [7, Def. 2.1].

We state our main result as follows.

THEOREM 2.1. *Consider an odd kernel $K_\epsilon(t)$, (2.3), normalized so that (2.4) holds, and satisfying the admissibility requirement (2.5). Then the kernel $K_\epsilon(t)$ satisfies the concentration property (2.2) for all piecewise smooth $f$'s, and the following error estimate holds:*

$$(2.7) \qquad |K_\epsilon * f(x) - [f](x)| \leq Const \cdot \epsilon.$$

---

[1]Here and below we use $BV[a,b]$ to denote the space of functions with *bounded variation*, endowed with the usual seminorm $||\phi||_{BV[a,b]} := \int_a^b |\phi'| dx.$

*Proof.* Using the fact that $K_\epsilon(t)$ is odd, we have

$$K_\epsilon * f(x) = -\int_{t\geq 0} K_\epsilon(t)(f(x+t) - f(x-t))dt$$

$$= -\int_{t\geq 0} tK_\epsilon(t)\frac{f(x+t) - f(x-t) - [f](x)}{t}dt - [f](x)\int_{t\geq 0} K_\epsilon(t)dt.$$

Applying (2.4) yields

$$(2.8) \qquad K_\epsilon * f(x) - [f](x) = -\int_{t\geq 0} tK_\epsilon(t)F_x(t)dt + \mathcal{O}(\epsilon).$$

By our assumption in (2.1), $F_x(t)$ is $BV$ and it is therefore bounded. Consequently, in the particular case that the moment bound (2.6) holds, the first term on the right of (2.8) is of order $\mathcal{O}(\epsilon)$, yielding

$$|K_\epsilon * f(x) - [f](x)| \leq Const \cdot \int |tK_\epsilon(t)|dt + \mathcal{O}(\epsilon) = \mathcal{O}(\epsilon).$$

In the general case, $F_x(t)$ has bounded variation, and the admissibility requirement (2.5) implies that the first term on the right of (2.8) is of order $\mathcal{O}(\epsilon)$, and we conclude

$$|K_\epsilon * f(x) - [f](x)| \leq Const||F_x(t)||_{BV} \cdot \epsilon + \mathcal{O}(\epsilon) = \mathcal{O}(\epsilon). \qquad \square$$

## 2.2. Examples of concentration kernels.

### 2.2.1. Compactly supported kernels.
Our first example consists of concentration kernels which "concentrate" near the origin, so that (2.6) holds. We consider a standard mollifier, $\phi_\epsilon(t) := \frac{1}{\epsilon}\phi(\frac{t}{\epsilon})$, based on an even, compactly supported bump function, $\phi \in C_0^1(-1, 1)$ with $\phi(0) = 1$. We then set

$$(2.9) \qquad K_\epsilon(t) = \frac{1}{\epsilon}\phi'\left(\frac{t}{\epsilon}\right) \equiv \phi'_\epsilon(t).$$

Clearly, $K_\epsilon$ is an odd kernel satisfying the required normalization (2.4)

$$\int_{t\geq 0} K_\epsilon(t)dt = \frac{1}{\epsilon}\int_{t\geq 0}\phi'\left(\frac{t}{\epsilon}\right)dt = -\phi(0) = -1.$$

In addition, its first moment is of order

$$\int_{t\geq 0}|tK_\epsilon(t)|dt = \epsilon\int_{0\leq s\leq 1}|s| \cdot |\phi'(s)|ds = \mathcal{O}(\epsilon),$$

and hence (2.6) holds. Theorem 2.1 then implies the following.

COROLLARY 2.2. *Consider the odd kernel* $K_\epsilon(t) = \phi'_\epsilon(t)$, *based on even* $\phi \in C_0^1(-1, 1)$ *with* $\phi(0) = 1$. *Then* $K_\epsilon(t)$ *satisfies the concentration property* (2.2), *and the following error estimate holds:*

$$(2.10) \qquad \phi'_\epsilon * f(x) = [f](x) + \mathcal{O}(\epsilon).$$

**2.2.2. The conjugate Dirichlet kernel.** The conjugate Dirichlet kernel,

$$K_N(t) = -\frac{1}{\log N}\tilde{D}_N(t), \qquad \tilde{D}_N(t) := \sum_{k=1}^{N} \sin kt,$$

is an example of an oscillatory concentration kernel. Clearly, $K_N(t)$ is an odd kernel. Moreover, the normalization (2.4) holds with $\epsilon \sim \frac{1}{\log N}$,

$$\int_{t=0}^{\pi} K_N(t)dt = -\frac{2}{\log N} \sum_{\text{odd } k\text{'s}} \frac{1}{k} = -1 + \mathcal{O}(\epsilon), \quad \epsilon \sim \frac{1}{\log N}.$$

Finally, summing

$$\tilde{D}_N(t) = \sum_{k=1}^{N} \sin kt = \frac{\cos\frac{t}{2} - \cos(N+\frac{1}{2})t}{2\sin\frac{t}{2}},$$

we find that the first moment of $K_N = -\tilde{D}_N(t)/\log N$ does not exceed

$$\int_{t=0}^{\pi} |tK_N(t)|dt \leq Const \cdot \epsilon, \quad \epsilon = \frac{1}{\log N},$$

so that the requirement (2.6) is fulfilled.

Theorem 2.1 then yields the classical result regarding the concentration of conjugate partial sums, [2, sect. 42], [23, sect. II, Thm. 8.13],

$$(2.11) \qquad -\frac{1}{\log N}\tilde{D}_N f(x) = [f](x) + \mathcal{O}\left(\frac{1}{\log N}\right).$$

We note in passing that in the case of the Dirichlet conjugate kernel, $K_N(t)$ does not concentrate near the origin, but instead (2.6) is fulfilled thanks to its uniformly small amplitude of order $\mathcal{O}(1/\log N)$. The error, however, is only of logarithmic order; consult [7, sect. 2].

**2.2.3. Oscillatory kernels: General concentration factors.** To accelerate the unacceptable logarithmically slow rate of the Dirichlet conjugate kernel in (2.11), we consider a general form of odd concentration kernels

$$(2.12) \qquad K_N^{\sigma}(t) := -\sum_{k=1}^{N} \sigma\left(\frac{k}{N}\right) \sin kt,$$

based on *concentration factors* $\sigma(\frac{k}{N})$ which are yet to be determined. Clearly $K_N^{\sigma}(t)$ is odd. Next, for the normalization (2.4) we note that

$$\int_{t=0}^{\pi} K_N^{\sigma}(t)dt = -2 \sum_{k \text{ odd}} \frac{\sigma(\frac{k}{N})}{k} \sim -\int_0^1 \frac{\sigma(x)}{x}dx.$$

In fact, the above Riemann's sum amounts to the midpoint quadrature, so that for $\frac{\sigma(\xi)}{\xi} \in C^2[0,1]$, one has

$$(2.13) \qquad \int_{t=0}^{\pi} K_N^{\sigma}(t)dt = -2 \sum_{k \text{ odd}} \frac{\sigma(\frac{k}{N})}{k} = -\int_0^1 \frac{\sigma(\xi)}{\xi}d\xi + \mathcal{O}\left(\frac{1}{N^2}\right),$$

and thus (2.4) holds for normalized concentration factors $\sigma(\xi)$,

$$(2.14) \qquad\qquad \int_0^1 \frac{\sigma(\xi)}{\xi} d\xi = 1.$$

Consult [7] for further refinement concerning the assumed regularity of $\sigma(\cdot)$ (We note that $\sigma(\cdot)$ is rescaled here with an additional factor of $-\pi$ compared to [7].)

Finally, we address the admissibility requirement (2.5) (and in particular (2.6)). To this end, we proceed along the lines of [7, Assertion 3.3], utilizing the identity (abbreviating $\xi_k = \frac{k}{N}$)

$$
\begin{aligned}
2\sin(t/2) K_N^\sigma(t) &\equiv \sum_{k=1}^{N-2} (\sigma(\xi_k) - 2\sigma(\xi_{k+1}) + \sigma(\xi_{k+2})) \frac{\sin(k+1)t}{2\sin(t/2)} \\
&\quad - (\sigma(1) - \sigma(\xi_{N-1})) \frac{\sin Nt}{2\sin(t/2)} + (\sigma(\xi_2) - \sigma(\xi_1)) \frac{\sin t}{2\sin(t/2)} \\
&\quad + \sigma(\xi_1) \cos\frac{t}{2} - \sigma(1) \cos\left(N + \frac{1}{2}\right)t \\
&=: I_1(t) + I_2(t) + I_3(t) + I_4(t) - I_5(t).
\end{aligned}
$$

(2.15)

This leads to the corresponding decomposition of $K_N^\sigma(t)$

$$K_N^\sigma(t) = R_N^\sigma(t) - \frac{\sigma(1)}{2} \frac{\cos(N + \frac{1}{2})t}{\sin\frac{t}{2}}.$$

Here, $R_N^\sigma(t)$ consists of the first four terms on the right-hand side of (2.15), $\sum_{j=1}^4 I_j(t)/2\sin(t/2)$, and it is easily verified that each one of these terms has a small first moment satisfying (2.6) (and consequently, (2.5) holds); i.e.,

$$(2.16) \qquad \int_{t=0}^\pi |t R_N^\sigma(t)| dt \le Const \cdot ||\sigma||_{C^2[0,1]} \frac{\log N}{N}.$$

For example, using the standard bound $|\sin(kt)/2\sin(t/2)| \le \min\{k, 1/t\}$, the contribution corresponding to the first term, $I_1(t)$, does not exceed

$$\int_{t=0}^\pi \left| t \frac{I_1(t)}{2\sin(t/2)} \right| dt \le \frac{\max|\sigma''|}{N^2} \sum_{k=1}^{N-2} \left[ \int_{t=0}^{1/N} + \int_{t=1/N}^\pi \right] \left| \frac{\sin(kt)}{2\sin(t/2)} \right| dt = \mathcal{O}\left( \frac{\log N}{N} \right).$$

Similar estimates hold for the remaining contributions of $I_2, I_3$, and $I_4$. In particular, since $\sigma(\xi)/\xi$ is bounded, $|\sigma(1/N)| \le \mathcal{O}(1/N)$, and hence $\int |t I_4(t)/2\sin(t/2)| dt = \mathcal{O}(1/N)$.

Finally, the admissibility of the fifth term on the right of (2.15) is due to standard cancellation which guarantees that (2.5) holds:

$$(2.17) \qquad \left| \frac{\sigma(1)}{2} \int_{t=0}^\pi t \frac{\cos(N + \frac{1}{2})t}{\sin\frac{t}{2}} \phi(t) dt \right| \le Const \cdot \frac{\sigma(1)}{N} ||\phi||_{BV}.$$

It is in this context of spectral concentration kernels that admissibility requires the more intricate property of cancellation of oscillations. Summarizing (2.13), (2.16), and (2.17), we obtain as a corollary an improved version of the main result in [7, Thm. 3.1] regarding spectral edge detection using concentration kernels, $K_N^\sigma(t)$. In particular, since $K_N^\sigma(t)$ are $N$-degree trigonometric polynomials, one detects the edges of the piecewise smooth function $f(x)$ directly from its spectral projection $S_N(f) := \sum_{k=-N}^{N} \hat{f}_k e^{ikx}$,

$$K_N^\sigma * f \equiv K_N^\sigma * S_N(f) = i\pi \sum_{k=-N}^{N} sgn(k)\sigma\left(\frac{|k|}{N}\right)\hat{f}_k e^{ikx}.$$

COROLLARY 2.3. *Consider the odd concentration kernel (2.12)*

$$K_N^\sigma(t) = -\sum_{k=1}^{N} \sigma\left(\frac{k}{N}\right)\sin kt, \qquad \frac{\sigma(\xi)}{\xi} \in C^2[0,1].$$

*Assume that $\sigma(\cdot)$ is normalized so that (2.14) holds:*

$$\int_0^1 \frac{\sigma(\xi)}{\xi}d\xi = 1.$$

*Then $K_N^\sigma(t)$ admits the concentration property (2.2), and the following estimate holds:*

$$(2.18) \qquad |K_N^\sigma * S_N(f) - [f](x)| = Const \cdot \frac{\log N}{N}.$$

*Remark.* One can relax the regularity on the concentration factor $\sigma(\cdot)$ [7]. Corollary 2.3 is a generalization of [7, Thm. 3.1];[2] in particular, the error estimate (2.18) is valid throughout the interval, including at the location of the jump discontinuities.

Let us introduce few prototypical examples of concentration factors $\sigma(\cdot)$ for the detection of edges from spectral data. In this context we note that other detection methods of discontinuities in periodic spectral data can be found in the works of Eckhoff [5], [6] and of Mhaskar and Prestin, e.g., [15] and the references therein. We note that our results apply to the nonperiodic expansions as discussed in section 3.2.2 below.

1. *Trigonometric factors.* We consider concentration factors of the form $\sigma(\xi) = \sigma_\alpha(\xi) := \sin\alpha\xi/Si(\alpha)$ with the proper normalization $Si(\alpha) := \int_0^\alpha (\sin\eta/\eta)d\eta$. The edge detector introduced originally by Banerjee and Geer [1] corresponds to $\sigma_\pi(\xi)$; the general case is found in [7, sect. 3.2].

2. *Polynomial factors.* As a first example consider $\sigma(\xi) = \xi$. In this case, $K_N^x * f = \frac{\pi}{N}S_N(f)'$, and Corollary 2.3 recovers Fejér's result, [23, sect. III, Thm. 9.3], with the following error estimate:

$$(2.19) \qquad \left|\frac{\pi}{N}S_N(f)' - [f](x)\right| \le Const \cdot \frac{\log N}{N}.$$

This is the first member of a whole family of polynomial concentration factors, e.g., [7, sect. 3.4],

$$(2.20) \qquad \sigma(\xi) = \sigma^p(\xi) := p\xi^p,$$

which correlate to concentration kernels satisfying (2.3), (2.4), and (2.5). For odd $p$'s, $K_N^{\sigma^p} * f = (-1)^{[p/2]}pN^{-p}S_N^{(p)}(f)$; for even $p$'s, $K_N^{\sigma^p} * f =$ $(-1)^{p/2}pN^{-p}H * S_N^{(p)}(f)$, where $H(x) = i \sum sgn(k)e^{ikx}$. These edge detectors were introduced in [9] and were recently analyzed by Kvernadze in [13]. Corollary 2.3 yields

$$\left| \frac{i\pi}{N^p} \sum_{k=-N}^{N} sgn(k)|k|^p \hat{f}_k e^{ikx} - [f](x) \right| \leq Const \cdot \frac{\log N}{N}.$$

---

[2]We note the different rescaling here of $\sigma(\cdot)$ by an additional factor of $-\pi$, compared with the formulation in [7, Thm. 3.1].

The last error estimate is (essentially) first order. It is sharp. It was noted in [7, sect. 3.4], however, that $\sigma^p$'s with higher $p$'s lead to faster convergence rates at selected interior points, bounded away from the singularities of $f$. This leads us to the next example.

3. *Exponential factors.* Polynomial concentration factors (of odd degree) correspond to differentiation in physical space; trigonometric factors correspond to divided differences in the physical space—consult the original derivation in [1]. Our main result stated in Corollary 2.3 provides us with the framework of general concentration kernels which are not necessarily limited to a realization in the physical space. In particular, we seek concentration factors, $\sigma(\cdot)$, which vanish at $\xi = 0, 1$ to any prescribed order,

$$(2.21) \qquad \frac{d^j}{d\xi^j}\sigma(\xi)\Big|_{\xi=0} = \frac{d^j}{d\xi^j}\sigma(\xi)\Big|_{\xi=1} = 0, \qquad j = 0, 1, 2, \ldots, p-1.$$

The higher $p$ is, the more localized the corresponding concentration kernel, $K_N^\sigma(\cdot)$, becomes. Here is why.

Evaluating $K_N^\sigma(t)$ at the equidistant points $t_\ell = 2\pi\ell/N$,

$$K_N^\sigma\left(\frac{2\pi\ell}{N}\right) = -\sum_{k=1}^{N} \sigma\left(\frac{k}{N}\right)\sin\frac{2\pi k\ell}{N},$$

we observe that $K_N^\sigma(t_\ell)/N$ coincide with the $\ell$-discrete Fourier coefficient of $\sigma(\cdot)$; since $\sigma(\xi)$ and its first $p$-derivatives vanish at both ends, $\xi = 0, 1$, there is a rapid decay of its (discrete) Fourier coefficients, $|\hat{\sigma}_\ell| \leq Const \cdot \ell^{-p}$,

$$|K_N^\sigma(t_\ell)| \leq Const \cdot \|\sigma\|_{C^{p+1}[0,1]} \frac{1}{(Nt_\ell)^p}.$$

Thus, for $t$ away from the origin, $K_N^\sigma(t)$ is rapidly decaying for large enough $N$'s. Moreover, we claim an increasing number of moments of $K_N^\sigma(\cdot)$ vanish. To this end we consider the odd moments of $K_N^\sigma(\cdot)$ (its even moments vanish, of course). With $\xi_k = k/N$ we find

$$(2.22) \qquad \int_{t=-\pi}^{\pi} t^j \times \left( -\sum_{k=1}^{N} \sigma\left(\frac{k}{N}\right)\sin kt \right) dt$$

$$= \frac{-2}{N^j} \sum_{k=1}^{N} \frac{1}{N} \int_{\tau=0}^{\pi N} \tau^j \sigma(\xi_k)\sin(\xi_k\tau)d\tau$$

$$= \sim \frac{\pm 2}{N^j} \int_{\xi=0}^{1} \sigma(\xi) \cdot \frac{d^j}{d\xi^j}\int_{\tau=0}^{\pi N}\cos(\xi\tau)d\tau d\xi, \quad j \text{ odd}.$$

Integrate by parts; respectively, sum by parts the summation on the right of (2.22). Thanks to (2.21) the boundary terms vanish and we have

$$\int_{t=-\pi}^{\pi} t^j K_N^\sigma(t)dt \sim \frac{1}{N^j}\int_{\xi=0}^1 \frac{d^j}{d\xi^j}\sigma(\xi)\cdot\frac{\sin(\pi N\xi)}{\xi}d\xi$$

$$= \frac{1}{N^j}\int_{\xi=0}^1 \frac{d^{p-j}}{d\xi^{p-j}}\Big(\frac{1}{\xi}\frac{d^j}{d\xi^j}\sigma(\xi)\Big)\frac{\sin(\pi N\xi)}{(\pi N)^{p-j}}d\xi = \mathcal{O}\left(\frac{1}{N^p}\right),$$

$$j < p, \ p \text{ odd.}$$

As an example, we consider the exponential concentration factors

$$(2.23)\quad \sigma^{exp}(\xi) = Const\cdot\xi e^{\frac{1}{\alpha\xi(\xi-1)}}, \qquad Const = \int \exp\left(\frac{-1}{\alpha\eta(\eta-1)}\right)d\eta$$

normalized so that $\int_{\xi=0}^1 \sigma^{exp}(\xi)/\xi d\xi = 1$. Here, the $C_0^\infty[-1,1]$ concentration factor $\sigma^{exp}(\xi)$ vanishes exponentially at both ends; $\xi = 0, 1$ so that (2.21) holds for *all* $p$'s. Figures 2.3 and 2.4 confirm the improved localization of these exponential concentration factors.

4. *Band pass filter.* Bauer [3] has considered a family of what he termed as "band pass filter," $\eta(\frac{k}{N})$, supported in the range of middle frequencies, say, supp$\eta \subset [1/4, 3/4]$. We note in passing that these are special cases of $p$-order admissible concentration factors (see (2.21)), although the normalization used in [3, eq. (1.35)], $\int \eta(x)/x \sin(\pi x)dx = 1$, prevented the recovery of the amplitude of the jumps.

To demonstrate the detection of edges by the concentration factors outlined above, consider the following two examples of discontinuous $f$'s (defined on $[-\pi, \pi]$):

$$f_a(x) := -sgnx\cdot\cos\left(\frac{x}{2}(2+sgnx)\right), \qquad f_b(x) := \begin{cases} \cos(x - \frac{x}{2}sgn(|x| - \frac{\pi}{2})), & x < 0, \\ \cos(\frac{5}{2}x + xsgn(|x| - \frac{\pi}{2})), & x > 0. \end{cases}$$

In both cases, $f_a(x)$ and $f_b(x)$ are recovered from their Fourier coefficients using the Fourier partial sums $S_N[f](x)$, and we wish to recover their jump discontinuities

$$[f_a](x) = \begin{cases} -2, & x = 0, \\ 0, & \text{else,} \end{cases} \qquad [f_b](x) = \begin{cases} \pm\sqrt{2}, & x = \pm\frac{\pi}{2}, \\ 0, & \text{else.} \end{cases}$$

Figures 2.1 and 2.2 demonstrate the use of trigonometric and polynomial concentration factors for the detection of edges from Fourier spectral data.

As noted in [10, sect. 3.4], polynomial factors of higher degree yield improved results away from the jump discontinuities. Indeed, the corresponding concentration kernels, $K_N^{\sigma^p}(\cdot)$, have additional vanishing moments. In the limit, one arrives at exponential factors, $K_N^{\sigma^{exp}}(\cdot)$; Figures 2.3 and 2.4 demonstrate the edge detection of these factors in Fourier expansions $S_N f_a(x)$ and $S_N f_b(x)$. The improved localization is evident, due to the faster convergence rate in the smooth parts of $f$'s. In particular, the superiority of the exponential factors is illustrated in the figures on the left, when compared with the first-order accurate polynomial concentration factor, $\sigma^{p=1}(\xi) = \xi$. At the same time, Gibbs oscillations can be noticed in the vicinity of the jump discontinuities.

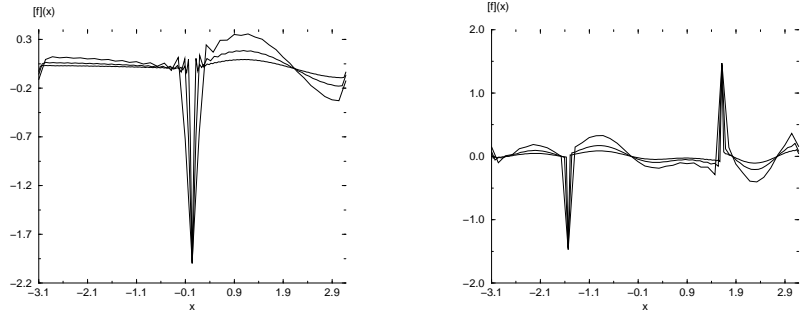

FIG. 2.1. *Trigonometric concentration factor* $\sigma(\xi) = \frac{\sin \xi}{Si(1)}$ *for (left)* $f_a(x)$, *where the exact jump value is* $[f](0) = -2$ *and (right)* $f_b(x)$, *where the exact jump values are* $[f](\pm\frac{\pi}{2}) = \pm\sqrt{2}$.

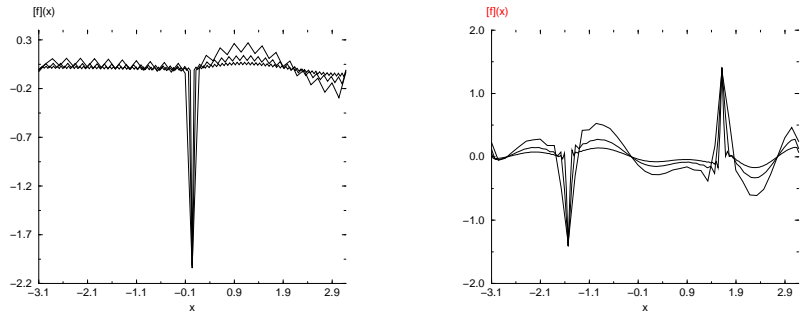

FIG. 2.2. *Jump value obtained by the polynomial concentration factor* $\sigma^{p=1}(\xi) = \xi$ *for (left)* $f_a(x)$ *and (right)* $f_b(x)$.

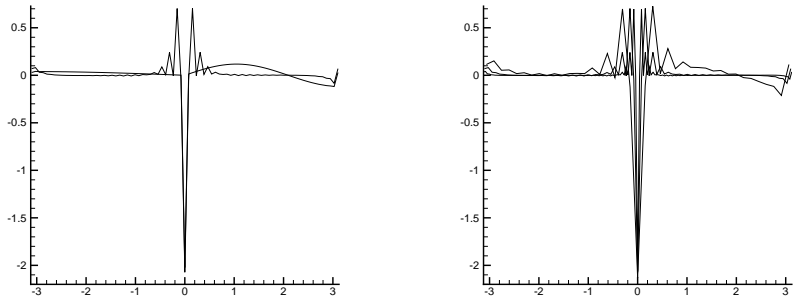

FIG. 2.3. *Edge detection using the exponential concentration factor* $\sigma(\xi) = 3\exp(\frac{1}{6\xi(\xi-1)})$ *(left)* $\sigma^{exp}$ *vs.* $\sigma^{p=1}$ *for* $S_{40}f_a(x)$ *and (right)* $\sigma^{exp}$ *for* $S_N f_a(x)$ *with* $N = 20, 40, 80$ *modes.*

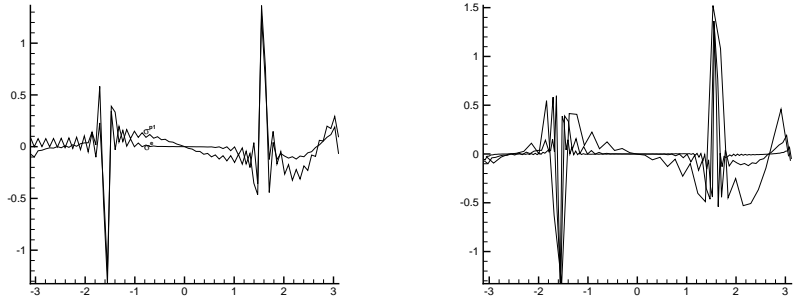

Fig. 2.4. *Edge detection using the exponential concentration factor* $\sigma(\xi) = 3\exp(\frac{1}{6\xi(\xi-1)})$ *(left)* $\sigma^{exp}$ *vs.* $\sigma^{p=1}$ *for* $S_{40}f_b(x)$ *and (right)* $\sigma^{exp}$ *for* $S_N f_b(x)$ *with* $N = 20, 40, 80$ *modes.*

**3. Edge detection in nonperiodic projections.** Consider a piecewise smooth $f(\cdot)$. To simplify our presentation, we assume $f$ experiences a single jump discontinuity at $x = c$. The localization property of the appropriate concentration kernel in the presence of a single jump applies to the case with finitely many jump discontinuities. We begin with an alternative derivation of our results for the periodic case.

**3.1. Revisiting the periodic case.** If a $2\pi$-periodic $f(\cdot)$ experiences a single jump, $[f](c)$, then it dictates the Fourier coefficients decay [4],[23],

$$\hat{f}_k = [f](c)\frac{e^{-ikc}}{2\pi ik} + \mathcal{O}\left(\frac{1}{k^2}\right).$$

To extract information about the location of the jump from the *phase* of the leading term, we examine the special concentration kernel, $K_N^\sigma$ with $\sigma(\xi) = \xi$, where $K_N^\xi * S_N(f) = \frac{\pi}{N}S_N(f)'$,

$$(3.1) \quad \frac{\pi}{N}S_N(f)' = \pi \sum_{k=-N}^{N} \frac{ik}{N}\left\{[f](c)\frac{e^{ik(x-c)}}{2\pi ik} + \mathcal{O}\left(\frac{1}{k^2}\right)\right\}$$

$$(3.2) \quad = [f](c)\frac{1}{2N}\sum\left\{1 + \mathcal{O}\left(\frac{1}{k}\right)\right\}e^{ik(x-c)} = [f](x) + \mathcal{O}\left(\frac{\log N}{N}\right).$$

Here we used the concentration property of the Dirichlet kernel localized at $x = c$:

$$\frac{1}{2N}\sum_{k=-N}^{N} e^{ik(x-c)} = \begin{cases} \mathcal{O}\left(\dfrac{1}{N|x-c|}\right), & x \neq c, \\[2mm] 1, & x = c. \end{cases}$$

The same property applies to the class of concentration factors, $\sigma(\xi) := \xi\mu(\xi)$, such that (2.14) holds, $\int_0^1 \mu(\xi)d\xi = 1$,

$$\frac{1}{2N}\sum_{k=-N}^{N} \mu\left(\frac{|k|}{N}\right)e^{ik(x-c)} = \begin{cases} \mathcal{O}\left(\dfrac{1}{N|x-c|}\right), & x \neq c, \\[2mm] 1, & x = c. \end{cases}$$

It then follows that the corresponding $K_N^\sigma$ in (2.12) is an admissible concentration kernel, so that $K_N^\sigma * f(x) \longrightarrow [f](x)$.

### 3.2. Nonperiodic expansions.

**3.2.1. General Jacobi expansions.** We begin with the Jacobi expansion of a piecewise smooth $f(\cdot)$,

$$(3.3) \qquad S_N(f) = \sum_{k=0}^N \hat{f}_k P_k(x), \qquad \hat{f}_k = \int_{-1}^1 f(x)\omega(x)P_k(x)dx.$$

Here $P_k(x)$ are the Jacobi polynomials—the eigenfunctions of the singular Sturm–Liouville problem

$$(3.4) \qquad ((1-x^2)\omega(x)P_k'(x))' = -\lambda_k \omega(x)P_k(x), \quad -1 \le x \le 1$$

with corresponding eigenvalues $\lambda_k = \lambda_k^{(\alpha)} = k(k+2\alpha+1)$. Different families of Jacobi polynomials are associated with different weight functions $\omega(x) \equiv \omega_\alpha(x) := (1-x^2)^\alpha$. To simplify the computations, we assume that the $P_k$'s are normalized so that $\|P_k(x)\|_\omega = 1$.

As in the periodic case, integration by parts (against (3.4)) shows that a single jump discontinuity, $[f](c)$, dictates the decay of the Jacobi coefficients,

$$(3.5)$$
$$\hat{f}_k = \frac{-1}{\lambda_k} \int_{-1}^1 f(x)((1-x^2)\omega(x)P_k'(x))'dx = [f](c)\frac{1}{\lambda_k}(1-c^2)\omega(c)P_k'(c) + \mathcal{O}\left(\frac{1}{\lambda_k^2}\right).$$

To extract information about the location of the jump, we consider the conjugate sum of the form

$$(3.6) \qquad \frac{\pi\sqrt{1-x^2}}{N}\sum_{k=1}^N \mu\left(\frac{k}{N}\right)\hat{f}_k P_k'(x)$$
$$= [f](c)\frac{\pi\sqrt{1-x^2}(1-c^2)\omega(c)}{N}\sum_{k=1}^N \mu\left(\frac{k}{N}\right)\left\{\frac{1}{\lambda_k} + \mathcal{O}\left(\frac{1}{\lambda_k^2}\right)\right\} \times P_k'(c)P_k'(x),$$

corresponding to concentration factors $\sigma(\xi) = \xi\mu(\xi)$. We shall focus our attention on the particular case $\mu(\xi) \equiv 1$,

$$(3.7)$$

$$\frac{\pi\sqrt{1-x^2}}{N}S_N(f)'(x) = \frac{\pi\sqrt{1-x^2}}{N}\sum_{k=1}^N \hat{f}_k P_k'(x)$$
$$= [f](c)\frac{\pi\sqrt{1-x^2}(1-c^2)\omega(c)}{N}\sum_{k=1}^N\left\{\frac{1}{\lambda_k} + \mathcal{O}\left(\frac{1}{\lambda_k^2}\right)\right\} \times P_k'(c)P_k'(x).$$

This is the nonperiodic analogue of the Fourier concentration kernel $K_N^\xi * f$ with the additional prefactor weight of $\sqrt{1-x^2}$.

We want to quantify the localization property of the last summation. To this end we note that if $\{P_k^{(\alpha)}(x)\}$ are the Jacobi polynomials w.r.t. the weight function $\omega_\alpha(x)$, then $\{P_k'(x)\}$ are the Jacobi polynomials w.r.t. the modified weight function $\omega_\beta(x) = (1-x^2)\omega_\alpha(x)$ with $\beta := \alpha + 1$: indeed, their $\omega_\beta$-orthogonality follows from integration by parts of (3.4) against $P_k^{(\alpha)}$. Thus,

$$P_k^{(\alpha)\prime}(x) = C_{k,\beta} P_{k-1}^{(\beta)}, \qquad \beta = \alpha + 1.$$

The coefficients $C_{k,\beta}$ are determined by normalization, where by using (3.4) once more we find

$$1 = \|P_k^{(\alpha)}\|_{\omega_\alpha}^2 = -\frac{1}{\lambda_k}\int_{-1}^{1}((1-x^2)\omega_\alpha(x)P_k^{(\alpha)\prime})'P_k^{(\alpha)}(x)dx = \frac{C_{k,\beta}^2}{\lambda_k}\|P_{k-1}^{(\beta)}\|_{\omega_\beta}^2,$$

and hence we set $C_{k,\beta} = \sqrt{\lambda_k}$, so that $\{P_{k-1}^{(\beta)}\}$ is the orthonormal family w.r.t. $\omega_\beta$ weight. Inserted into the leading term of (3.7), we end up with a Jacobi kernel associated with weight function $\omega_\beta(x) = (1-x^2)^\beta$,

$$\frac{\pi\sqrt{1-x^2}}{N}S_N(f)'(x) \sim [f](c)\frac{\pi\sqrt{1-x^2}\omega_\beta(c)}{N} \times \sum_{k=1}^{N}\frac{C_{k,\beta}^2}{\lambda_k}P_{k-1}^{(\beta)}(c)P_{k-1}^{(\beta)}(x)$$

$$= [f](c)\frac{\pi\sqrt{1-x^2}\omega_\beta(c)}{N} \times \sum_{k=1}^{N}P_{k-1}^{(\beta)}(c)P_{k-1}^{(\beta)}(x).$$

We rewrite this as

$$(3.8) \qquad \frac{\pi\sqrt{1-x^2}}{N}S_N(f)'(x) \sim [f](c)\frac{\pi\sqrt{1-x^2}\omega_\beta(c)}{N} \times K_N(c,x).$$

By virtue of the Christoffel–Darboux formula, e.g., [19, Thm. 3.2.2], the kernel $K_N(c,x)$ is given by

$$(3.9) \qquad K_N(c,x) = \frac{k_{N-1}}{k_N}\frac{P_N^{(\beta)}(x)P_{N-1}^{(\beta)}(c) - P_N^{(\beta)}(c)P_{N-1}^{(\beta)}(x)}{x-c}, \qquad \frac{k_{N-1}}{k_N} \sim \frac{1}{2},$$

and it remains to quantify the concentration property of $K_N(c,x)$. To this end we use the asymptotic behavior of $P_N^{(\beta)}$ which is stated as[3]

$$(3.10) \quad P_N^{(\beta)}(x) \sim \min\left(\sqrt{\frac{2}{\pi\omega_{\beta+1/2}(x)}}, Const \cdot N^{\beta+1/2}\right) = \sqrt{\frac{2}{\pi\omega_{\beta+1/2}(x_N)}},$$

where $x_N := sgn(x) \cdot \min\{|x|, 1 - Const/N\}$ denotes the separation between the interior and boundary regions. Using this to upper bound $K_N(c,x)$ in (3.9), we find

$$(3.11) \qquad |K_N(c,x)| \leq \frac{1}{\pi\sqrt{\omega_{\beta+1/2}(c_N)\omega_{\beta+1/2}(x_N)}} \times \frac{1}{|x-c|}, \quad x \neq c.$$

---

[3]The first term on the right of (3.10) follows from the classical asymptotic formula, e.g., [19, Thm. 12.1.4], which tells us the behavior of the $L_\omega^2$-normalized $P_N^{(\beta)}(x)$ at the interior $x = \cos\theta$

$$\sqrt{\omega_{\beta+1/2}(x)}P_N^{(\beta)}(x) \sim \sqrt{2/\pi}\cos\{N\theta + \gamma(\theta)\}, \quad x = \cos\theta.$$

The second term on the right of (3.10) reflects the fact that as $x$ approaches the $\pm1$-boundaries, the $L_\omega^2$-normalized $P_N^{(\beta)}(x)$ approaches its maximal value, e.g., [19, eqs. (4.7.3), (4.7.15)]

$$P_N^{(\beta)}(1) = \binom{N+2\lambda-1}{N}\sqrt{\frac{\Gamma(N+1)(N+\lambda)}{\pi 2^{1-2\lambda}\Gamma(N+2\lambda)}}\Gamma(\lambda) \sim \sqrt{\frac{\Gamma(N+2\lambda)(N+\lambda)}{N!}} \sim N^\lambda, \quad \lambda = \beta + 1/2.$$

The upper bound on the right is in fact the leading term in the asymptotics of $K_N(c, x)$ for large $N$'s as long as $-1 < x \neq c < 1$, e.g., [19, sect. 13.4]. Similarly, the behavior at $x = c$ yields

$$(3.12) \qquad K_N(c, c) \sim \frac{N + Const}{\pi \omega_{\beta+1/2}(c)}, \quad -1 < c < 1.$$

The desired concentration property now follows, similar to the localization of the periodic Dirichlet kernel $D_N(x - c)/N$ in (3.3). We restrict our attention to *interior jumps*, $-1 < c < 1$, so that for $N$ large enough, $c_N = c$, and (3.11), (3.12), and (3.8) yield

$$(3.13)$$

$$\frac{\pi\sqrt{1-x^2}\,\omega_\beta(c)}{N} \times K_N(c, x) \sim \begin{cases} \mathcal{O}\left( \dfrac{\sqrt{\omega_{\alpha+1/2}(c)}}{\sqrt{\omega_{\alpha+1/2}(x_N)}} \times \dfrac{1}{N|x-c|} \right) \sim \dfrac{N^{\alpha+1/2}}{N}, & x \neq c, \\ \\ 1, & x = c. \end{cases}$$

We summarize by stating the following.

COROLLARY 3.1. *Let $S_N(f)$ denote the truncated Jacobi expansion (3.3) of a piecewise smooth $f$, associated with a weight function $\omega_\alpha = (1-x^2)^\alpha$, $-1 < \alpha \leq 0$. Then $\pi\sqrt{1-x^2}S_N(f)'(x)/N$ admits the concentration property*

$$(3.14) \qquad \left| \frac{\pi\sqrt{1-x^2}}{N} S_N(f)'(x) - [f](x) \right| \leq Const \cdot \frac{\log N}{N(1-x^2)^{\alpha/2+1/4}},$$

$$-1 + Const \cdot \frac{1}{N^2} < x < 1 - Const \cdot \frac{1}{N^2}.$$

It is instructive to examine the above discussion for the special case of Chebyshev expansion corresponding to $\alpha = -\frac{1}{2}$,

$$S_N(f)(x) = \sum_{k=0}^{N} {}' \hat{f}_k T_k(x), \quad \hat{f}_k = \frac{2}{\pi} \int_{x=-1}^{1} \frac{f(x)T_k(x)}{\sqrt{1-x^2}} dx.$$

(Observe that except for Chebyshev expansion, the concentration bound (3.13) deteriorates as we approach the boundaries, depending on whether $|x| \sim 1$ for $\alpha > -\frac{1}{2}$ or $|c| \sim 1$ for $\alpha < -\frac{1}{2}$.) The conjugate sum corresponding to (3.8) reads

$$\frac{\pi\sqrt{1-x^2}}{N} S_N(f)'(x) \sim [f](c)\frac{\pi}{N} \sum_{k=1}^{N} {}' \frac{2}{\pi k^2} \sqrt{1-c^2}\, T_k'(c) \sqrt{1-x^2}\, T_k'(x).$$

In this case, we can sum the corresponding Chebyshev kernel: setting $x = \cos\theta$ and $-1 < c = \cos\eta < 1$ we find

$$(3.15)$$

$$\frac{\pi\sqrt{1-x^2}}{N} S_N(f)'(x) \sim [f](c)\frac{2}{N} \sum_{k=0}^{N} {}' \sin(k\eta)\sin(k\theta)$$

$$= [f](c)\frac{1}{N} \Big[ D_N(\theta - \eta) - D_N(\theta + \eta) \Big] = \begin{cases} \mathcal{O}(\frac{1}{N}), & x \neq c, \\ \\ [f](c), & x = c. \end{cases}$$

**3.2.2. Chebyshev expansion.** Our discussion above on edge detection in the nonperiodic expansions is based on expansion of the Jacobi coefficients to their leading order in (3.5). More precise information is obtained using the general framework introduced in the main theorem, Theorem 2.1.

COROLLARY 3.2. *Let $f(\cdot)$ be a piecewise smooth function with Chebyshev expansion $S_N f(x) \sim \sum' \hat{f}(k) T_k(x)$. Consider the concentration factors, $\sigma(\xi) := \xi\mu(\xi)$, with $\mu(\cdot)$ normalized so that*

$$\int_0^1 \mu(\xi)d\xi = 1, \qquad \mu(\xi) \in C^2[0,1].$$

*Then $K_N^\sigma(t) * f(cos\theta)$ admits the concentration property (2.2), and the following estimate holds:*

$$(3.16) \qquad \left| \frac{\pi\sqrt{1-x^2}}{N} \sum_{k=1}^N \mu\left(\frac{k}{N}\right) \hat{f}(k) T_k'(x) - [f](x) \right| = Const \cdot \frac{\log N}{N}.$$

*Proof.* With a piecewise smooth $f(x)$ defined over the interval $[-1,1]$ we utilize the usual Chebyshev transformation $x = \cos\theta$, $0 \le \theta \le \pi$. We consider the *even* extension $f(\cos\theta)$, $-\pi \le \theta \le \pi$. Using Theorem 2.1 along the lines of Corollary 2.3, we find that the odd concentration kernel, $K_N^\sigma(t)$, recovers the jumps of $f(\cos\theta)$, i.e.,

$$\left| -\sum_{k=1}^N \frac{k}{N}\mu\left(\frac{k}{N}\right)\sin(kt) * f(\cos\theta) - [f](\cos\theta) \right| \le Const \cdot \frac{\log N}{N}.$$

With $T_k'(x) = -k\sin k\theta/\sqrt{1-x^2}$, a straightforward computation shows the sum on the left equals

$$-\sum_{k=1}^N \sigma\left(\frac{k}{N}\right)\sin(kt) * f(\cos\theta) = -\pi\sum_{k=1}^N{}'\frac{k}{N}\mu\left(\frac{k}{N}\right)\hat{f}(k)\sin(k\theta)$$

$$= \frac{\pi\sqrt{1-x^2}}{N}\sum_{k=1}^N \mu\left(\frac{k}{N}\right)\hat{f}(k)T_k'(x)$$

and the result follows. □

We turn to numerical examples. The following tables summarize our results for the edge detection in Legendre expansion, corresponding to $\alpha = 0$, and in Chebyshev expansion, corresponding to $\alpha = -1/2$. Scaled to the unit interval $[-1,1]$, we consider $f_a(\frac{x}{\pi})$ and $f_b(\frac{x}{\pi})$. The results confirm the linear convergence rate stated in Corollary 3.1, both away from the jumps—consult Tables 3.1 and 3.2—as well as at the jump itself (Table 3.3).

We note that the critical threshold must be very high for $N = 20, 40$ to eliminate the artificial jumps. This indicates that 40 nodes are not enough to resolve the jumps of $f_b(x)$ in either the Chebyshev or Legendre case.

It is clear from Tables 3.1 and 3.2 that convergence is nonuniform at the boundaries. We have observed in our numerical experiments that the edge detector, $\frac{\pi\sqrt{1-x^2}}{N}S_N[f]'(x)$, experiences larger oscillations near the boundaries which do affect the linear convergence rate there. In this context we note the dependence of the error bounds on the smoothness of $f(\sqrt{1-x^2})$.

The first-order convergence is reconfirmed, in Table 3.4 below, when measuring the $L^1$-error away from the jumps discontinuities (and *up* to the boundaries).

TABLE 3.1

*Pointwise error estimate $|\pi\sqrt{1-x^2}/NS_N(f_a)'(x/\pi) - [f_a](x/\pi)|$ away from the jump discontinuity at $x = 0$.*

| $N$ | Legendre expansion | | | | Chebyshev expansion | | | |
|---|---|---|---|---|---|---|---|---|
| | $x=-.998$ | $x=-.9805$ | $x=.5$ | $x=.75$ | $x=-.998$ | $x=-.9805$ | $x=.5$ | $x=.75$ |
| 40 | .192 | .112 | .135 | .123 | 3.6E-2 | 7.1E-2 | .121 | .147 |
| 80 | 7.7E-2 | 3.8E-2 | .16 | 1.3E-2 | 2.1E-2 | 5.9E-3 | .145 | 2.4E-2 |
| 160 | 7.6E-3 | 4.0E-3 | 3.1E-2 | 5.6E-3 | 5.3E-3 | 8.7E-3 | 3.3E-2 | 1.6E-2 |

TABLE 3.2

*Pointwise error estimate $|\pi\sqrt{1-x^2}/NS_N(f_b)'(x/\pi) - [f_b](x/\pi)|$ away from the jump discontinuities at $x = \pm\frac{1}{2}$.*

| $N$ | Legendre expansion | | | | Chebyshev expansion | | | |
|---|---|---|---|---|---|---|---|---|
| | $x=-.998$ | $x=-.9805$ | $x=.5$ | $x=.75$ | $x=-.998$ | $x=-.9805$ | $x=.5$ | $x=.75$ |
| 40 | 8.0E-2 | 7.4E-2 | 5.4E-3 | .45 | 1.3E-02 | 1.5E-02 | 9.8E-02 | .56 |
| 80 | .15 | 1.55E-2 | 7.0E-4 | .18 | 1.1E-02 | 1.3E-02 | 4.5E-02 | .26 |
| 160 | 1.6E-3 | 7.2E-3 | 3.9E-4 | .12 | 1.8E-03 | 4.6E-03 | 2.3E-02 | .12 |

TABLE 3.3

*Pointwise error estimate $|\pi\sqrt{1-x^2}/NS_N(f)'(c) - [f](c)|$ at the point(s) of discontinuity, $x = c$.*

| | Legendre | | | Chebyshev | | |
|---|---|---|---|---|---|---|
| $N$ | $f_a(0)$ | $f_b(-\pi/2)$ | $f_b(\pi/2)$ | $f_a(0)$ | $f_b(-\pi/2)$ | $f_b(\pi/2)$ |
| 40 | 1.9E-02 | 8.5E-02 | .187 | 1E-2 | 5.3E-2 | 6.2 E-2 |
| 80 | 1.4E-03 | 3.9E-02 | .11 | 1E-2 | 4.7E-2 | 2.5E-2 |
| 160 | 1.3E-02 | 1.7E-02 | 5.5E-02 | 8.7E-3 | 1.1E-2 | 2.2E-2 |

TABLE 3.4

*$L^1[-1,1] - \{c\}$ error estimate $|\pi\sqrt{1-x^2}/NS_N(f)'(x) - [f](x)|$ away from discontinuities.*

| $N$ | Legendre | | Chebyshev | |
|---|---|---|---|---|
| | $f=f_a$ | $f=f_b$ | $f=f_a$ | $f=f_b$ |
| 40 | .16 | .25 | .156 | .27 |
| 80 | 8.3E-2 | .13 | 8.4E-2 | .12 |
| 160 | 4.4E-2 | 6.8E-2 | 4.6E-2 | 6.9E-2 |

**4. Nonlinear enhancement.** The detection of edges in Theorem 2.1 is based on *separation of scales*. Thus, consider for example a piecewise smooth $f$ with finitely many jump discontinuities at $x = c_1, c_2, \ldots, c_n$. If $K_\epsilon$ is an admissible concentration kernel, then $|K_\epsilon * f(x)| \ll 1$ for $x$ away from these jumps, whereas at $x = c_j$, $K_\epsilon * f(x)_{x=c_j} \sim O(1)$,

$$(4.1) \qquad K_\epsilon * f(x) = \begin{cases} \mathcal{O}(\epsilon), & x \neq c_1, c_2, \ldots, \\ [f](c_j), & x = c_j. \end{cases}$$

The last statement refers to the asymptotic behavior of the concentration kernel as a function of the small parameter $\epsilon \downarrow 0$. In this section we outline a new, nonlinear enhancement procedure, which is easily implemented to "pinpoint" finitely many edges in piecewise smooth $f$'s.

To this end we *enhance* the separated scales in (4.1) by considering

$$(4.2) \qquad \epsilon^{-\frac{p}{2}}(K_\epsilon * f(x))^p = \begin{cases} \sim Const \cdot \epsilon^{\frac{p}{2}}, & x \neq c_1, c_2, \ldots, \\ ([f](c_j))^p \cdot \epsilon^{-\frac{p}{2}}, & x = c_1, c_2, \ldots. \end{cases}$$

By increasing the exponent $p > 1$, we enhance the separation between the vanishing scale at the points of smoothness (of order $\mathcal{O}(\epsilon^{\frac{p}{2}})$) and the growing scale at the jumps (of order $\mathcal{O}(\epsilon^{-\frac{p}{2}})$).

Next, one must introduce a *critical threshold* which will eliminate all the unacceptable jumps. Only those edges with amplitudes larger than the critical threshold, $[f](x) > J_{crit}^{1/p}\sqrt{\epsilon}$, will be detected. Thus $J = J_{crit}$ is a measure which defines the small scale in our computation of edge detection. We note that $J = J_{crit}$ is data dependent and is typically related to the variation of the smooth part of $f$.

Given this critical threshold, we form our *enhanced concentration kernel* $K_{\epsilon,J}[f]$

$$(4.3) \qquad K_{\epsilon,J}[f](x) = \begin{cases} K_\epsilon * f(x) & \text{if } \epsilon^{-\frac{p}{2}}|K_\epsilon * f(x)|^p > J_{crit}, \\ 0 & \text{otherwise.} \end{cases}$$

Clearly, with $p$ large enough, one ends up with a sharp edge detector where $K_{\epsilon,J}[f](x) = 0$ at all but $\mathcal{O}(\epsilon)$ neighborhoods of the jumps $x = c_1, c_2, \ldots$ In practical applications, a moderate enhancement exponent, $p \leq 5$, will suffice. We consider two examples.

1. *The quadratic filter.* Consider the peaked concentration kernel (2.9) $K_\epsilon(t) = \phi'_\epsilon(t)$. Then, with $p = 2$, one finds the so-called *quadratic filter* [12], [22], where

$$(4.4) \qquad (K_\epsilon * f(x))^2 = (\phi'_\epsilon * f(x))^2 \to [f]^2(x).$$

2. *Enhanced spectral concentration kernels.* We apply the procedure of nonlinear enhancement in conjunction with spectral concentration kernels $K_\epsilon = K_N^\sigma(t) := -\sum_{k=1}^N \sigma(\frac{k}{N}) \sin kt$ by considering the corresponding *enhanced spectral concentration kernel*

$$K_{N,J}^\sigma = \begin{cases} K_N^\sigma * S_N f(x) & \text{if } \epsilon^{-\frac{p}{2}}|K_N^\sigma * S_N(x)|^p > J_{crit}, \\ 0 & \text{otherwise.} \end{cases}$$

The enhanced spectral concentration kernel depends on four ingredients which are at our disposal:
- The number of modes, $N$;
- The enhancement exponent, $p$;
- The critical threshold, $J$;
- The concentration factor, $\sigma(\xi)$.

Figures 4.1 and 4.2 demonstrate the enhancement procedure to the spectral detection of edges depicted earlier in the corresponding Figures 2.2 and 2.1.

We conclude with nonperiodic examples. In Figure 4.3 we show the detection of a single edge in $f_a(x/\pi)$ from its Legendre expansion, $\sum (\hat{f}_a)_k P_k(x)$. The detection in Chebyshev expansion is shown in Figure 4.4 for $f_b(x/\pi)$. In both cases we used an enhancement factor $p = 2$ and a critical threshold $J = 5$.

**5. Concluding remarks.** Accurate reconstruction of piecewise smooth functions from their spectral projections is plausible only when the location (and amplitude) of the underlying jump discontinuities are known; consult [1], [7], [5], [6], [10], [16], [21], [11], and references therein.

Theorem 2.1, and its Corollaries 2.3 and 3.1, provide the general framework for the detection of edges from spectral data, in both periodic and nonperiodic cases. The detection is based on admissible concentration kernels which include as particular cases classical examples of Fejér as well as additional examples in recent literature

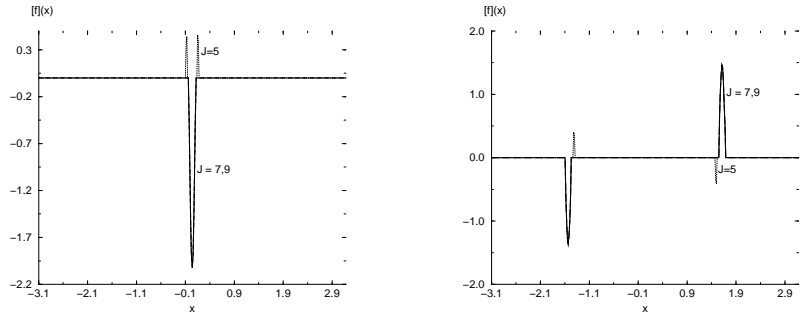

FIG. 4.1. *Jump value obtained by applying the polynomial concentration factor* $\sigma(\xi) = \xi$ *with* $N = 40$ *with enhancement exponent* $p = 2$ *for* (a) $f_a(x)$, *where the exact jump value is* $[f](0) = -2$ *and* (b) $f_b(x)$, *where the exact jump values are* $[f](\pm\frac{\pi}{2}) = \pm\sqrt{2}$.

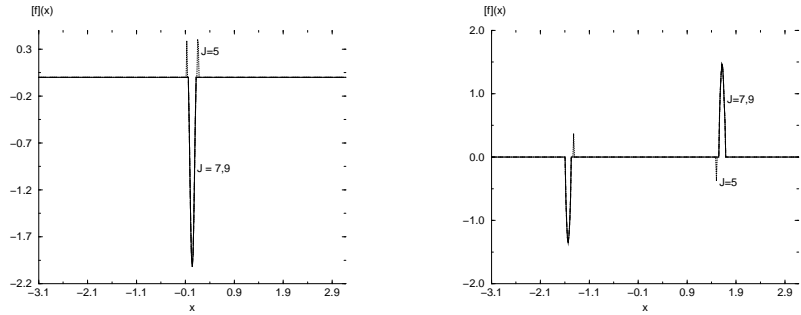

FIG. 4.2. *Jump value obtained by applying the trigonometric concentration factor* $\sigma_1(\xi) = \frac{\sin\xi}{Si(1)}$ *with* $N = 40$ *modes and enhancement exponent* $p = 2$ *for* (left) $f_a(x)$, *where the exact jump value is* $[f](0) = -2$ *and* (right) $f_b(x)$, *where the exact jump values are* $[f](\pm\frac{\pi}{2}) = \pm\sqrt{2}$.

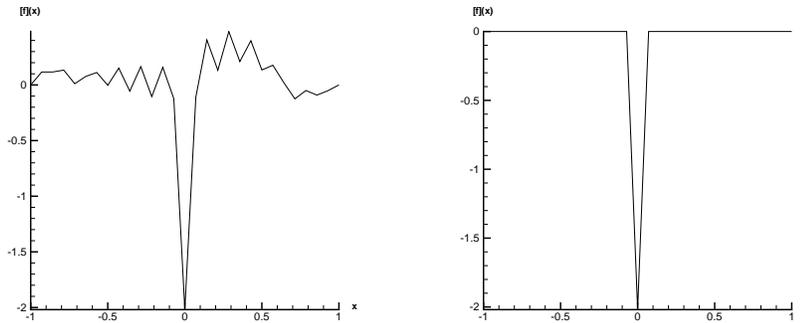

FIG. 4.3. *Detection of edges in Legendre expansion of* $f_a(x/\pi)$ *with exact jump value is* $[f_a](0) = -2$ (left) *before and* (right) *after enhancement with* $p = 2$ *and* $J = 5$.

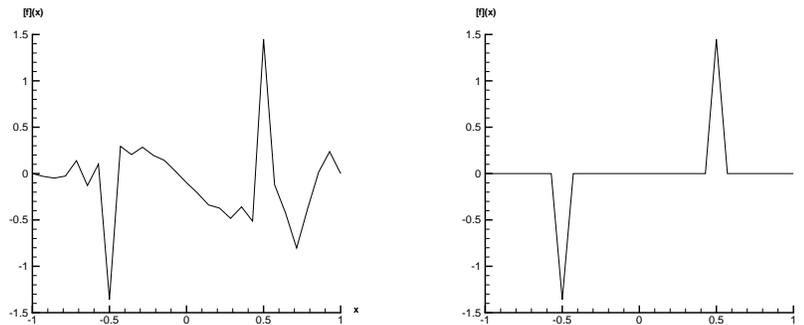

FIG. 4.4. *Detection of edges in Chebyshev expansion of $f_b(x/\pi)$ with exact jump value is $[f_b](\pm\frac{1}{2}) = \pm\sqrt{2}$ (left) before and (right) after enhancement with $p = 2$ and $J = 5$.*

[1], [9], [13]. In particular, we introduce here a new family of exponential concentration kernels, (2.23), with a superior convergence rate away from the edges. A linear convergence rate is observed near the detected edges. We also introduce a nonlinear enhancement (4.3) procedure which enables one to "pinpoint" edges with amplitude larger than a critical threshold.

Recently, the edge detection and enhancement method was applied to nonlinear conservation laws [20], [8] as a postprocessing tool to improve the overall convergence rate of the spectral viscosity solution. Since the edge detection occurs only at the postprocessing stage, very little cost is added to the procedure yet the results are dramatically improved. Future applications, in both one- and several-space dimensions, will also include image processing, where edge detection is needed to denoise the contamination by the $\mathcal{O}(1)$-Gibbs oscillations in the neighborhoods of the undetected edges.